%


\input mssymb


%

%

\newcount\skewfactor
\def\mathunderaccent#1#2 {\let\theaccent#1\skewfactor#2
\mathpalette\putaccentunder}
\def\putaccentunder#1#2{\oalign{$#1#2$\crcr\hidewidth
\vbox to.2ex{\hbox{$#1\skew\skewfactor\theaccent{}$}\vss}\hidewidth}}


\def\eqdef{\buildrel \rm def \over =}

\def\phi{\varphi}

\def\sat{\models}
\def\su{\subseteq}
\def\a{\alpha}
\def\b{\beta}

\def\d{\delta}
\def\l{\lambda}

\def\om{\omega}

\def\lng{\langle}
\def\rng{\rangle}
\def\ov{\overline}

\def\cont{{2^{\aleph_0}}}



\def\inv{{\rm Inv}}

\def\otp{{\rm otp}\,}

\def\id{{\rm id}}
\def\Inv{{\rm INV}}

\def\im{{\rm Im}}


\def\endproof#1{\hfill  
{\parfillskip0pt$\smiley_{\hbox{{#1}}}$\par\medbreak}}

\def\imply{\Rightarrow}
\def\iff{\Leftrightarrow}
\def\proof{\smallbreak\noindent{\sl Proof}: }



\newbox\noforkbox \newdimen\forklinewidth
\forklinewidth=0.3pt   

\setbox0\hbox{$\textstyle\bigcup$}
\setbox1\hbox to \wd0{\hfil\vrule width \forklinewidth depth \dp0
			height \ht0 \hfil}
\wd1=0 cm
\setbox\noforkbox\hbox{\box1\box0\relax}
\def\unionstick{\mathop{\copy\noforkbox}\limits}
\def\nonfork#1#2_#3{#1\unionstick_{\textstyle #3}#2}
\def\nonforkin#1#2_#3^#4{#1\unionstick_{\textstyle #3}^{\textstyle  
#4}#2}

\setbox0\hbox{$\textstyle\bigcup$}
\setbox1\hbox to \wd0{\hfil$\nmid$\hfil}
\setbox2\hbox to \wd0{\hfil\vrule height \ht0 depth \dp0 width
				\forklinewidth\hfil}
\wd1=0cm
\wd2=0cm
\newbox\doesforkbox
\setbox\doesforkbox\hbox{\box1\box0\relax}
\def\nunionstick{\mathop{\copy\doesforkbox}\limits}

\def\fork#1#2_#3{#1\nunionstick_{\textstyle #3}#2}
\def\forkin#1#2_#3^#4{#1\nunionstick_{\textstyle #3}^{\textstyle  
#4}#2}

\font\circle=lcircle10

\setbox0=\hbox{~~~~~}
\setbox1=\hbox to \wd0{\hfill$\scriptstyle\smile$\hfill} 
\setbox2=\hbox to \wd0{\hfill$\cdot\,\cdot$\hfill} 

\setbox3=\hbox to \wd0{\hfill\hskip4.8pt\circle i\hskip-4.8pt\hfill}  


\wd1=0cm
\wd2=0cm
\wd3=0cm
\wd4=0cm

\newbox\smilebox
\setbox\smilebox \hbox {\lower 0.4ex\box1
		 \raise 0.3ex\box2
		 \raise 0.5ex\box3
		\box4
		\box0{}}
\def\smiley{\leavevmode\copy\smilebox}

\headline={\tenrm  
\number\folio\hfill\jobname\hfill\number\day.\number\month.\number\ye 
ar}


\outer\long\def\ignore#1\endignore{}

\newcount\itemno
\def\itm{\advance\itemno1 \item{(\number\itemno)}}
\def\ritm{\advance\itemno1 \item{)\number\itemno(}}
\def\startitm{\itemno=0 }
\def\aitm{\advance\itemno1 

\item{(\letter\itemno)}}
^^L

\def\letter#1{\ifcase#1 \or a\or b\or c\or d\or e\or f\or g\or h\or
i\or j\or k\or l\or m\or n\or o\or p\or q\or r\or s\or t\or u\or v\or
w\or x\or y\or z\else\toomanyconditions\fi}
\def\raitm{\advance\itemno1 \item{)\rletter\itemno(}}
\def\rletter#1{\ifcase#1\or `\or a\or b\or c\or d\or e\or f\or g\or
h\or i\or k\or k\or l\or n\or n\or q\or r\or t\or v
\else\toomanyconditions\fi}


\newcount\secno
\newcount\theono

\catcode`@=11
\newwrite\mgfile

\openin\mgfile \jobname.mg
\ifeof\mgfile \message{No file \jobname.mg}
	\else\closein\mgfile\relax\input \jobname.mg\fi
\relax
\openout\mgfile=\jobname.mg

\newif\ifproofmode
\proofmodetrue            

\def\@nofirst#1{}

\def\neusection{\advance\secno by 1\relax \theono=0\relax}
\def\neuchap{\secno=0\relax\theono=0\relax}

\neuchap

\def\labelit#1{\global\advance\theono by 1%
             \global\edef#1{%
             \number\secno.\number\theono}%
             \write\mgfile{\@definition{#1}}%
}

^^L




\def\ppro#1#2:{%
\labelit{#1}%
\smallbreak\noindent%
\@markit{#1}%
{\bf\ignorespaces {#2}:}}





\def\@definition#1{\string\def\string#1{#1}
\expandafter\@nofirst\string\%
(\the\pageno)}

\def\@markit#1{
\ifproofmode\llap{{ \expandafter\@nofirst\string#1\ }}\fi%
{\bf #1\ }
}

\def\h@markit#1{
\ifproofmode\edef\nxt{\string#1\ }%
{\tenrm\beginL\nxt\endL}
\fi%
{{\bf\beginL #1\endL}}
}
 ^^L
\def\labelcomment#1{\write\mgfile{\expandafter
		\@nofirst\string\%---#1}}

\catcode`@=12

\def\refrence#1#2:{\write\mgfile{\def\noexpand#1{#2}}%
\areference{#2}}
\def\areference#1{\medskip\item{[#1]} \ignorespaces}

\newcount\referencescount
\def\numrefrence#1{\advance\referencescount1
\edef#1{\number\referencescount}%
\write\mgfile{\def\noexpand#1{#1}}%
\areference{#1}}

\def\numericalreferences{\let\refrence\numrefrence}

\newcount\scratchregister
\def\simplepro{\scratchregister\theono\advance\scratchregister1 

\edef\scratchmacro{\number\secno.\number\scratchregister}%
\expandafter\ppro\csname\scratchmacro\endcsname}


\font\teneuf=eufm10
\font\seveneuf=eufm7
\font\fiveeuf=eufm5
\newfam\euffam
\textfont\euffam=\teneuf
\scriptfont\euffam=\seveneuf
\scriptscriptfont\euffam=\fiveeuf

\overfullrule=0pt
\magnification=1200
\baselineskip=21pt
\centerline{\bf Universal Abelian Groups}
\centerline{by Menachem Kojman and Saharon
Shelah\footnote{*}{Partially supported by the United
States--Israel Binational science foundation. Publication number  
455}}
\medskip
\proofmodefalse
\headline{}
\centerline{ABSTRACT}
{ We examine the existence of universal elements in classes of
infinite abelian groups. The main method is using group invariants
which are defined relative to club guessing sequences. We prove, for
example:

\noindent
{\bf Theorem}: For $n\ge 2$, there is a purely universal
separable $p$-group in $\aleph_n$ if, and only if, $\cont\le
\aleph_n$.                                             

}
\bigskip
{\bf \S\number\secno\ Introduction}

In this paper ``group'' will always mean ``infinite abelian group'', 

and ``cardinal'' and ``cardinality''  always refer to  infinite
cardinals
and infinite cardinalities.

Given a class of groups $K$ and a cardinal $\l$ we call a group $G\in
K$ {\bf universal for $K$ in $\lambda$} if $|G|=\l$ and every $H\in  
K$
with $|H|\le\l$ is isomorphic to a subgroup of $G$.  The objective of
this paper is to examine the existence of universal groups in various
well-known classes of infinite abelian groups.  We also investigate
the existence of {\bf purely universal} groups for $K$ in $\l$,
namely groups $G\in K$ with $|G|=\l$ such that every
$H\in K$ with $|H|\le \l$ is isomorphic to a {\bf pure} subgroup of
$G$.

The main set theoretic tool we use is a {\it club guessing sequence}.
This is a prediction principle which has enough power to control
properties of an infinite object which are defined by looking at all
possible enumerations of the object. Unlike the diamond and the
square, two combinatorial principles which are already accepted as
useful for the theory of infinite abelian groups, club guessing
sequences are proved to exist in ZFC.  Therefore using club guessing
sequences does not require any additional axioms beyond the usual
axioms of ZFC. Club guessing sequences are particularly useful in
proving theorems from {\bf negations} of CH and GCH.
\bigbreak

The paper is organized as follows: in section 1 we define 

group invariants relative to  club guessing sequences, and show
that the invariants are  monotone in pure embeddings. In Section 2 we
construct various groups with prescribed demands on their invariants.
In section 3 these ideas are used to investigate the existence of
universal groups for classes of torsion groups and classes of torsion
free groups. It appears that cardinal arithmetic decides the question
of existence of a universal group in  many cardinals. For example:
there is a purely universal separable $p$-group in $\aleph_n$ iff
$\aleph_n\ge \cont$ for all $n\ge 2$ (for $n=1$ only the ``if'' part
holds).

This paper follows two other papers by the same authors, [KjSh 409]
and [KjSh 447], in which the existence of universal linear orders,
boolean algebras, and models of unstable and  stable unsuperstable  
first order
theories were examined using the same method.

 All the abelian group theory one needs here, and more, is found in
[Fu],  whose system of notation we adopt.  An acquaintance with
ordinals and cardinals is necessary, as well as familiarity with
stationary sets and the closed unbounded filter.  Knowledge of  
chapter
II in [EM] is more than enough.

Before getting on, we first observe that in every infinite  
cardinality
there are
universal groups which are divisible:

\ppro \triv Theorem:
In every cardinality there is a universal group, universal $p$-group
(for every
prime $p$), universal torsion group and universal torsion-free group.

\proof These are, respectively, the direct sum of $\l$ copies of the
rational group $Q$ together with $\l$ copies of $Z(p^\infty)$ for
every prime $p$; the direct sum of $\l$ copies of $Z(p^\infty)$; the
direct sum of $\l$ copies of $Z(p^\infty)$ for every prime $p$; and
the direct sum of $\l$ copis of $Q$. The universality of these  
groups,
each for its respective class, follows from the structure theorem for
divisible groups and the fact that every group ($p$-group, torsion
group, torsion-free group) is embeddable in a divisible group
($p$-group, torsion group, torsion-free group) of the same  
cardinality
([Fu] I, 23 and 24])\endproof\triv

\neusection 

\noindent{\bf \S\number\secno\ The invariant of a group relative to
the ideal id$(\ov C)$}

A fixed assumption in this section is that $\l$ is a regular
uncountable cardinal. We assume the reader is familiar with the basic
properties of closed unbounded sets of $\l$, and with the definition
and basic properties of stationary sets.

\ppro \congDef Definition: For a group $G$, $nG\eqdef \{ng:g\in G\}$.
Two elements $g,h\in G$ are {\bf $n$-congruent} if $g-h\in nG$. If
$g,h$ are $n$-congruent, we also say that $h$ is an {\bf
$n$-congruent} of $g$.

\ppro \pureDef Definition: 

\item{(1)}([Fuch, p.113]) Let $G$ be  group. A subgroup $H\su G$ is a
{\bf 

pure} subgroup, denoted by $H\su_{pr} G$,  if for all natural $n$,
$nH=nG\cap H$. 

\item{(2)} An embedding of groups $h:H\to G$ is a {\bf pure}
embedding if its image $h(H)$ is a pure subgroup of $G$.

\ppro \representation Definition: 

 Suppose that $\l$ is a regular uncountable cardinal and that
$G$ is a  group of cardinality $\l$. A sequence $\ov G=\lng
G_\a:\a<\l\rng$ is called a {\bf $\l$-filtration } of $G$  iff for  
all $\a$ 

\startitm
\itm $G_\a\su G_{\a+1}$
\itm $G_\a$ is of cardinality smaller than $\l$
\itm if $\a$ is limit, then
$G_\a=\bigcup_{\b<\a}G_\b$
\itm $G=\bigcup\limits_{\a<\l} G_\a$. 


$G$,
obtained

 Suppose $\ov G=\lng G_\a:\a<\l\rng$ is a given filtration
of a group $G$.  Suppose $c\su \l$ is a set of ordinals, and
the increasing enumeration of $c$ is $\lng \a_i:i<i(*)\rng$. Let
$g\in G$ be an element. We define a way in which  $g$ chooses a  
subset
of $c$:

\ppro \invDef Definition:
$\inv_{\ov G}(g,c)=\{\a_i\in c: g\in \bigcup_n ((G_{\a_i+1}+nG) -
(G_{\a_i} + nG))\}$

We call $\inv_{\ov G}(g,c)$ {\bf the
invariant of the element $g$ relative to the $\l$-filtration $\ov G$
and the set of indices $c$}.

Worded otherwise,  $\inv_{\ov G} (g,c)$, is the subset of those  
indices
$\a_i$ such
that by increasing the group $G_{\a_i}$ to the larger group
$G_{\a_{i+1}}$, 

an $n$-congruent  for $g$  is
introduced for some $n$.

As the definition of the invariant depends on a $\l$-filtration, one
may think that the invariant does not deserve its name.
Indeed, given a group $G$ equipped with two respective
$\l$-filtrations $\ov G$ and $\ov G'$, it 

is not necessarily true that for $g\in G$
$$\inv_{\ov G}(g,c)=\inv_{\ov G'}(g,c)\leqno{(1)}$$

The solution to this problem is working with a club guessing sequence
$\ov C=\lng c_\d:\d\in S\rng$ and the ideal $\id(\ov C)$ associated  
to
it. The idea is as follows: for any pair of $\l$-filtrations $\ov G$
and $\ov G'$ a group $G$ there is a club  $E\su
\l$ such that for every $\a\in E$, $G_\a=G'_\a$.  So if we chose our
set $c$ in the definition of invariant to consist only of such
``good'' $\a$-s, namely if $c\su E$, then it does not matter
according to which $\l$-filtration we work.  But we cannot choose a
set $c$ which is a subset of every club $E$ resulting from some
pair of $\l$-filtrations. What we {\bf can} do is find a sequence of
$c$-s with the property that for every club $E\su
\l$, stationarily many of them are subsets of $E$. Thus we will be  
able
to define an invariant that is  independent of a
particular choice of a $\l$-filtration.
Here is the precise formulation of this:

\ppro \clubGuessDef Definition: A sequence $\lng c_\d:\d\in S\rng$,
where $S\su \l$ is a stationary set, $c_\d\su \d$ and $\d=\sup c_\d$
for every $\d$, is called a {\bf club guessing sequence} if for every
club $E\su \l$ the set $\{\d\in S:c_\d\su E\}$ is a stationary 

subset of $\l$.

The theorems asserting the existence of club guessing sequences will  
be
quoted later. A club guessing sequence $\ov C=\lng c_\d:\d\in S\rng$
gives rise to an ideal $\id(\ov C)$ over $\l$ --- the {\bf guessing
ideal}.

\ppro \guessIdeal Definition: Suppose that $\ov C=\lng c_\d :\d\in
S\rng$ is a club gussing sequence. We define a proper ideal $\id(\ov
C)$ as follows:
\item{} $A\in \id(\ov C)\iff A\su \l \;\&\; \exists E\su \l$, $E$  
club,
$\&\;\forall \d\in E\cap S, c_\d\not\su E$.

So a set of $\d$-s is small if there is a club $E$ which it fails to
guess stationarily often, namely there is no $\d\in A$ such that
$c_\d\su E$.

\ppro \ideal Lemma:
If $\ov C$ is a club guessing sequence as above, then $\id(\ov C)$ is  
a
proper, $\l$ complete ideal.

\proof That $\id(\ov C)$ is proper means that it does not contains
every subset of $\l$. Indeed, $S\notin \id(\ov C)$, as it guesses
every club. That $\id(\ov C)$ is downward closed is immediate from  
the
definition. Suppose, finally, that $A_i, i<(*)$, $i(*)<\l$ are sets  
in
the ideal. We show that their union $A\eqdef\bigcup\limits
_{i<i(*)}A_i$ is in the ideal. Pick a club $E_i$ for every $i<i(*)$  
so
that $\d\in A_i \imply c_\d\not\su E_i$. The set $E=\bigcap\limits
_{i<i(*)}E_i$ is a club.  Suppose that $d\in A$. Then there here is
some $i<i(*)$ such that $\d\in A_i$.  Therefore $c_\d\not\su E_i$.  
But
$E\su E_i$, so necessarily $c_\d\not\su E$. Thus, $A\in\id(\ov
C)$\endproof\ideal

We adopt the phrase ``for almost every $\d$ in $S$'', by which we  
mean  ``for all
$\d\in S$ except for a set in $\id (\ov C)$''.

\ppro \wellDef Lemma:  Suppose that $\ov C = \lng c_\d:\d\in S\rng$  
is
a
club guessing sequence on  $S\su \l$. Suppose that $\ov G$ and $\ov
G'$ are two $\l$-filtrations of a group $G$ of cardinality $\l$. Then
for
almost every $\d\in S$, (1) holds for every $g\in G$.

\proof: The set of $\a<\l$ for which $G_\a=G'_a$ is a club. Let us
denote it by $E$. If for some $\d$,  $c_\d\su E$ holds, then
for every $g\in G$ it is true that $\inv_{\ov G}(g,c_\d)=\inv_{\ov
G'}(g,c_\d)$. But as $\ov C$ is a club guessing sequence, by
definition, for almost every $\d$, $c_\d\su E$.\endproof{\wellDef}

We define now the desired  group invariant.

\ppro \secInvDef Definition: Suppose that $\ov C$ is a club guessing
sequence and that $\ov G$ is a $\l$-filtration of a group $G$ of
cardinality $\l$. Let
\startitm
\itm $P_\d(\ov G,\ov C)=\{\inv_{\ov G}(g, c_\d):g\in G\}$
\itm $\Inv(G,\ov C)=[\lng P_\d(\ov G, \ov C):\d\in S\rng]_{\id(\ov  
C)}$

 The second item should read ``the
equivalence class of the sequence of $P_\d$ modulo the ideal $\id(\ov
C)$'', where two sequences are equivalent modulo an ideal if the set  
of
coordinate in which the sequences differ is in the ideal.

\ppro\veryWellDef Lemma: The definition of $\Inv(G,\ov G)$ does not
depend on the choice  $\l$-filtration.

\proof Suppose that $\ov G, \ov G'$ are two $\l$-filtrations. By the
regularity of $\l$, there exists a club
$E$ such that for every $\a\in E$, $G_\a=G'_\a$. Therefore for every
$\d$ such that $c_\d\su E$ and every $g\in G$, $\inv_{\ov
G}(g,c_\d)=\inv_{\ov G'}(g,c_\d)$. This means that for every $\d$  
such
that $c_\d\su E$, $P_\d(\ov G,\ov C)=P_\d(\ov G',\ov C)$. But for  
almost
all $\d$ it is true that $c_\d\su E$, therefore the sequqnces $\lng
P_\d(\ov G,\ov C):\d\in S\rng$  and   $\lng P_\d(\ov G',\ov C):\d\in
S\rng$ are equivalent modulo $\id(\ov C)$.\endproof\veryWellDef

We remark at this point that the definition just made depends on the
existence of a club guessing sequence!  Strangely enough, we can  
prove
the existence of club guessing sequences for regular uncountable
cardinals $\l$ for all such cardinals {\bf except} $\aleph_1$.

Let us now quote the relevant theorems which assert the existence of
club guessing sequences:

\ppro\clubGuessExist Theorem: If $\mu$ and $\l$ are cardinals,
$\mu^+<\l$ and $\l$ is regular, then there is a club guessing  
sequence
$\lng c_\d:\d\in S\rng$ such that the order type of each $c_\d$ is
$\ge\mu$.

\proof In [Sh-e, new VI\S2]= [Sh-e, old III\S7].

We procced to show that  $\Inv$ is preserved, in a
way, under pure embeddings.

\ppro \embedding Lemma: Suppose that $H$ and $G$ are groups of
cardinality $\l$ and that $\ov H$ and $\ov G$ are $\l$-filtrations.
Suppose that $\ov C$ is a club guessing sequence on $S\su \l$.  If
$h:H\to G$ is a pure embedding, then for almost every $\d\in S$
$P_\d(\ov
H,\ov C)\su P_\d(\ov G,\ov C)$.

\proof
Suppose for simplicity that $H\su_{pr} G$, namely that the embedding
is the identity function. The set $E_1\eqdef \{\a:H\cap G_\a=H_\a\}$  
is
a
club.   Define
for every natural number $n$  a function $f_n(y)$ on $G$ as follows:

$$f_n(y)=\cases{{\rm some} \quad x\in\{x:
x\in H\,\&\,  (x+y)\in nG\} & if
$\{x:x\in H\,\&\, (x+y)\in nG\}\not=\emptyset$ \cr
 & \cr
0 & otherwise}
$$

There is a club $E_2$ such that $G_\a$ is closed under $f_n$ for
all $n$ for every $\a\in E_2$. $E=E_1\cap E_2$ is a club.

\ppro \good Claim: Suppose that $h\in H$ and that $\a\in E$. Then $h$
has an $n$-congruent 

in $G_\a$ (in the sense of $G$) iff $h$ has an $n$-congruent in  
$H_\a$
(in the sense of $H$).

\proof
One direction is trivial. Suppose, then,   that there is an  
$n$-congruent
$g\in G_\a$.  Let $h'=f_n(g)$. By the definition of $f_n$, $h'+g\in
nG_\a$; also $h-g\in nG$. Therefore $h-h'\in nG$. As $H\su_{pr} G$,  
$h-
h'\in nH$, and therefore $h'$ is an $n$-congruent of $h$ in the sense
of $H$. \endproof{\good}

The proof of the Lemma follows now readily: For almost every $\d\in  
S$
it is true $c_\d\su E$. Therefore for every such $\d$, every $h\in H$
and every $n$, $h$ has an $n$-congruent in $H_\a$ iff $h$ has an
$n$-congruent in $G_\a$. Therefore $\inv _{\ov H}(h,c_\d)\su\inv_{\ov
G}(h,c_\d)$. \endproof{\embedding}

\bigbreak
\bigbreak
\noindent
\neusection
{\bf \S\number\secno\ Constructing groups with prescribed INV }

 In
this section we construct several groups with prescribed demands
on their $\Inv$. These groups are used in the next section to  show
that in certain cardinals universal groups do not exist. The method
in all constructions is  attaching to a simply defined group
points from a topological completion of the group. 

\bigbreak
\noindent
{\bf a. Constructions of $p$-groups}

\ppro \firstPco Theorem: If $\l$ is a regular uncountable cardinal,
$\ov C=\lng c_\d: \d\in S\rng$ is a club guessing sequence 

and $A_\d\su c_\d$ is a given
set of order type $\om$, then there is a separable $p$-group $G$  of
cardinality $\l$  and $\l$-filtration $\ov G$ such that $A_\d\in
P_\d(\ov G,\ov C)$ for every $\d\in S$.

\ppro \Premark Remark: This implies by Lemma \embedding\ that for
every separable $p$-group $G'$ of cardinality $\l$ which purely
extends $G$  and a $\l$-filtration $\ov G'$, for
almost every $\d$, $A_\d\in P_\d(\ov G',\ov C)$

\proof For every $n$,
let $B_n= \bigoplus\limits_{\eta\in \,{}^n\l}A_\eta$ where $A_\eta$  
is a copy
of $Z_{p^n}$ with generator $a_\eta$. Let  
$G^0=\bigoplus\limits_{n}B_n$,
and let $G^1$ be the torsion completion of $G^0$. $G^1$ may be
identified with all sequences $(x_1,x_2,\cdots)$ where $x_n\in B_n$
and such that there is a (finite) bound to $\{o(x_n)\}_n$.  For  
details see
[Fuchs II,14--21]. The group we seek lies between $G^0$ and $G^1$,  
and
is a pure subgroup of $G^1$.

Let us make a simple observation:

\startitm
\itm If $x=(x_1,x_2,\cdots)$ and $y=(y_1,y_2,\cdots)$ belong to $G^1$
and $x-y\in p^nG^1$, then $x_i=y_i$ for all $i\le n$.

\proof Let $z_i=x_i-y_i$. $z_i\in B_i$. As $B_i\cap p^nG^1=0$ for  
$i\le
n$, we are done.

For every $\d\in S$ let $\lng \a^\d_n:n\in \om\rng$ be the increasing
enumeration of $A_\d$. Denote by $\eta^\d_n$
the sequence $\lng \a^\d_1,\cdots,\a^\d_n\rng$. Let $b^0_\d\in G^1$  
be
$(x^\d_1,x^\d_2,\cdots)$ where $x^\d_n=p^{n-1}a_{\eta^\d_n}$. So  
$x_n$ is of
order $p$ and height $n-1$. Consequently, $b^0_\d$ is of order $p$. 

  Let us denote 

$${x^\d_k\over {p^n}}\eqdef\cases
{p^{k-n-1}a_{\eta^\d_k} & if $n<k$\cr
&\cr
0& otherwise}
$$ and also let ${0\over {p^n}}\eqdef 0$.
Let
$b^n_\d=(\cdots,{{x^\d_k}\over {p^n}}\cdots)$. Let
$G$ be the subgroup of $G^1$ generated by $G^0$ together with
$\{b^n_\d:\d\in S,\;n<\om\}$. Having defined $G$, let us specify a
$\l$-filtration $\ov G$. For every $i<\l$ let $G_i$ be $\lng
\{a_\eta:\eta\in {}^{<\om}i\}\cup \{b^n_\d:\d\in  
S,\;\d<i,\;n<\om\}\rng$.


\ppro \auxClaim Claim: $\inv_{\ov G}(b^0_\d,c_\d)=A_\d$.

\proof  We should show that the set of indices $i$ with the
property that in $G_{i+1}$ some congruent of $b^0_\d$ appears  
coincides
with $A_\d$.  Suppose first, then, that $i=\a_n$ for some $n$.
$(x^\d_1,\cdots,x^\d_n,0,0,\cdots)$ is clearly a $p^n$-congruent of
$b^\d_\d$, as $b^0_\d-(x^\d_1,\cdots,x^\d_n)=p^nb^n_\d$.
Conversely, suppose that $i<\a_n$ and suppose to the contrary that  
there
is some
$y=(y_1,y_2,\cdots)\in G_i$ such that $b^0_\d-y\in P^nG$. Then  
$y_i=x_i$
for all $i\le n$ by (1). But $y\in G_i$ implies that $y_n\in G_i$ ---  
a
contradiction to $i<\a_n$. \endproof{\auxClaim,\firstPco}

\noindent
{\bf b. Constructions of torsion-free groups}

We start by constructing a torsion-free homogeneous group of a given
type ${\bf t}=(\infty,\cdots,\infty,0,\infty\cdots)$.  We recall that
a characteristic $\chi(g)$ of an element $g\in G$ is the sequence
$(k_1,k_2,\cdots)$ where $k_l$ is the $p$-height of $g$ for the  
$l$-th
prime. A  height can be $\infty$. A type ${\bf t}$ is an equivalence
class of characteristics modulo the equivalence relation of having
only a finite
difference in a finite number of coordinates. A homogeneous group is  
a
group in which all elements have the same type. We call a type ${\bf
t}$ a $p$-type if ${\bf t}=(\infty,\cdots,\infty,0,\infty\cdots)$
where the only coordinate in which there is $0$ is the number of $p$
in the list of primes.

\ppro \typeCo Theorem:
For every uncountable and regular cardinal $\l$, a club guessing
sequence $\ov C=\lng c_\d:\d\in S\rng$ and given sets $A_\d\in c_\d$,
each $A_\d$ of order type $\om$, there is a homogeneous group $G$ of
cardinality $\l$ with $p$-type ${\bf t}$
and a $\l$-filtration $\ov G$ such that for every $\d\in S$, $A_\d\in
P_\d(\ov G,\ov C)$.

\ppro \typeRemark Remark: This means that for every pure extension
$G'$ of $G$, for almost every $\d\in S$, $A_\d\in P_\d(\ov G',\ov  
C)$.

\proof
This proof resembles the proof of Theorem \firstPco. Let
$G^0=\bigoplus\limits_\l Q_p$ (where $Q_p$ is the group or rationals
with denominators prime to $p$). We index the isomorphic copies of  
$Q_p$
by $\eta\in {}^n\l$ and fix $a_\eta$, an element $a_\eta$ of  
characteristic 

$(\infty,\cdots,\infty,0,\infty\cdots)$ in the $\eta$-th copy of
$Q_p$. Let $G^1$ be the completion of $G^0$ in the $p$-adic topology. 

Let $\lng \eta^\d(n):n<\om\rng$ be the increasing enumeration of  
$A_\d$,
and let $\eta^\d_n=\lng \eta^\d(0),\cdots,\eta^\d(n-1)\rng$. Let
$b_{\d,n}=\sum\limits_k p^{k-n} a_{\eta_k}$. The rest is as in the  
proof of
Theorem \firstPco.\endproof{\typeCo}

\neusection
\noindent
{\bf\S\number\secno\ The main results}

\noindent
{\bf a. The Universality Spectrum of Torsion groups}

There is universal torsion group in $\l$ iff there is a universal
$p$-group in $\l$ for every prime $p$. We therefore may focus on
$p$-groups alone. There is a universal divisible $p$-group in $\l$,
the group $\bigoplus\limits_\l Z(p^{\infty})$, therefore the first
interesting question to ask in torsion groups is whether there is  a
universal {\bf reduced} $p$-group. Here the answer is ``no'':

\ppro \rankNeg Theorem: If $\l$ is an infinite cardinal (not
necessarily regular, not necessarily uncountable) then there is no
universal reduced $p$-group in $\l$.

\proof  There are $p$-groups of cardinality $\l$ of Ulm length  
$\sigma$
for every ordinal $\sigma<\l^+$. As $u(A)\le u(B)$ whenever $A\su B$,
and $u(A)<\l^+$, for every group of cardinality $\l$, no $p$-group of
cardinality $\l$ can be universal.\endproof\rankNeg

We put a further restriction on the class of $p$-groups, by demanding
that the Ulm length of a group be at most  
$\omega$.\footnote{$^1$}{One
can make finer distictions here by considering the class of all
$p$-groups of Ulm length which is bounded by an ordinal $\sigma$. But
we do not do this here.} We restrict ourselves then to the class of
separable $p$-groups. On this class see [Fu] vol II, chapter XI.

\noindent
{\bf b.  Universal separable $p$-groups}

We investigate the universality spectrum of the class or separable
$p$-groups.

\ppro\sepPPos Theorem: If $\l=\l^{\aleph_0}$ then there is a purely
universal separable $p$-group in $\l$.

\proof
Let $B=\bigoplus B_n$ where $B_n=\bigoplus\limits_{\l}Z_{p^n}$. The
torsion completion of $B$, denoted by $G$, is of cardinality
$|B|^{\aleph_0}=\l^{\aleph_0}=\l$ and is puely universal in $\l$. To
see this let $A$ be any separable $p$-group of cardinality $\l$, and
let $B_A$ be its basic subgroup. $B_A$ is purely embeddable in $B$,  
and
this gives rise to a pure embedding of $A$ into $G$.\endproof\sepPPos

We see then, that for every $n$ such that $\aleph_n\ge \cont$ there  
is
a purely universal separable $p$-group in $\aleph_n$. As CH implies
that in every $\aleph_n^{\aleph_0}=\aleph_n$ for all $n$, it follows
by \sepPPos\ that there is a purely universal separable $p$-group in
every $\aleph_n$. It is not uncommon that CH decides questions in
algebra. It is much less common, though, that a negation of CH or of
GCH does the same. The following theorem uses a negation of GCH as  
one
of its hypotheses.

\ppro \sepPNeg Theorem: $\l$ is regular and there is some $\mu$ such
that
$\mu^+<\l<\mu^{\aleph_0}$ then there is no purely universal separable
$p$-group in $\l$.

\proof
By $\mu^+<\l$ and Theorem \clubGuessExist, we may pick some club
guessing sequence $\ov C=\lng c_\d:\d\in S\rng$ where $S$ is a  
stationary
set of $\l$ and $\otp c_\d\ge\mu$. Suppose $G$ is a given separable
$p$-group. We will show that $G$ is not universal by presenting a
separable $p$-group $H$ of cardinality $\l$ which is not embeddable  
in
$G$. We choose a $\l$-filtration $\ov G$ of $G$ and observe that
$|P_\d(\ov G,\ov C)|\le \l$ for every $\d$. As $\mu^{\aleph_0}>\l$,
there is some $A_\d\su c_\d$,  of order type $\om$, which
dose not belong to $P_\d(\ov G,\ov C)$. By Theorem \firstPco, there  
is a
group $H$ of cardinality $\l$ such that for every embedding  
$\phi:H\to
G$, for almost every $\d$, $A_\d\in P_\d(\ov G,\ov C)$. This can hold
only
emptily, that is, if there are no such embeddings, because $A_\d$ was
chosen such that $A_\d\notin P_\d(\ov G,\ov C)$\endproof\sepPNeg

\ppro\chCor Corollary: For $n\ge 2$, there is a purely universal
separable $p$-group in $\aleph_n$ if, and only if, $\cont\le
\aleph_n$.

\proof: If $\aleph_n\ge 2^{\aleph_0}$ then
$\aleph_n^{\aleph_0}=\aleph_n$ and by Theorem \sepPPos\ there is a
purely universal separable $p$-group in $\aleph_n$. Conversely, if
$n\ge 2$ and $\aleph_n<\cont$, then by \sepPNeg\ there is no purely
universal separable $p$-group in $\aleph_n$\endproof{\chCor}

\noindent
{\bf c.  The Universality Spectrum of Torsion-Free Group}

We may restrict discussion in this Section to reduced
torsion-free groups.
We proceed to show first that in regular $\l$ which satisfy
$\l=\l^{\aleph_0}$ there is a universal reduced torsion-free group.
The proof of the next theorem is an isolated point in the paper with
respect to the technique, because it employs model theoretic
notions (first order theory, elementary embedding and saturated
model). These are available in every standard textbook on model
theory, like [CK].

\ppro\torFreePos Theorem: if $\l=\l^{\aleph_0}\ge 2^{\aleph_0}$, then
there
is a universal reduced torsion-free group in $\l$.

\proof Let $T$ be a complete first order theory of torsion
free-groups. It is enough to find a reduced group $G_T$ of  
cardinality
$\l$ such that $G_T\sat T$ and for every $H\sat T$, $H$ is embedded  
in
$G_T$; for if we have such a $G$ for every $T$, the group
$\bigoplus\limits_{T}G_T$ is of cardinality $\l$ (there are only
$2^{\aleph_0}$ complete first order theories), and is evidently
universal.

Let, then, $G'_T$ be a saturated model of $T$ of cardinality $\l$.  
Let
$D$ be its maximal division subgroup, and let $G_T\eqdef G'_T/D$.  
$G_T$
is
isomorphic to the direct summand of $D$, and is therefore  
torsion-free
and reduced. Suppose that $H\sat T$ is reduced (and, clearly,
torsion-free). There is an elementary embedding $f:H\to G'_T$.

\ppro\divClaim Claim: $\im f\cap D=0$

\proof Suppose $0\not=a\in H$ and $f(a)\in D$. As $f$ is elementary,
$a$ is divisible in $H$ by every integer $n$. As $H$ is torsion-free,
the set of all divisors of $a$ generates  a divisible subgroup of  
$H$,
contrary to $H$ being reduced.

We conclude, therefore, that $\hat f$ defined by $\hat f(a)=f(a)+D$
is an embedding of $H$ into $G_T$.\endproof{\torFreePos}

Next we show that below the continuum
there is no purely-universal reduced torsion-free group. The reason  
for
this
is trivial: there are $\cont$ types (over the empty set) in this  
class.
Therefore we do not need the club guessing machinery, and gain an  
extra
case -- the case where $\l=\aleph_1$ --- in comparison to Corollary  
\chCor.

\ppro \firstTorFreeNeg Theorem: If $\l<\cont$ then there is no
purely-universal reduced torsion-free group in cardinality $\l$. In
fact, for
every reduced torsion-free group $G$ of cardinality $\l$ there is a
rank-1 group $R$ which is not purely embeddable in $G$.

\proof 

As $\l<\cont$, there is a characteristic $(k_1,k_2\cdots)$, with all
$k_i$ finite, which is not equal to  $\chi_G(g)$ for every $g\in G$
(The definitions of {\bf characteristic} and {\bf type} are from
[Fu], II, 85). Let $R$ be
a rank-1 group such that $\chi_R(1)=(k_1,k_2\cdots)$. As pure
embeddings preserve the characteristic, $R$ is not purely embeddable
in $G$.\endproof\firstTorFreeNeg

We look now at a lasse of torsion-free groups which do not have many
types above the empty set. This is the class of homogeneous groups (a
group is {\bf homogeneous} if all non zero elements in the group have
the same type. See [Fu] II p.109).  Here we are able again to prove
that there is no purely universal group in the class in cardinality
$\l<\cont$ {\bf if} $\l>\aleph_1$. However, rather than using
types over the empty set, we are using invariants.

\ppro \homNeg Theorem: If $\l$ is a regular cardinal,
$\mu^+<\l<\mu^{\aleph_0}$ for some $\mu$, and ${\bf t}$ is a given
$p$-type, then there is no purely universal torsion free group in  
$\l$. Even more, for every torsion free group of cardinality $\l$  
there is a homogeneous
for the class of homogeneous groups whose type is ${\bf t}$.

\proof Let $G$ be any homogeneous group with type ${\bf t}$, and fix
some $\l$-filtration  $\ov G$. Let $\ov C$ be a club guessing  
sequence,
and for every $\d\in S$ let $A_\d$ be such that $A_\d\notin P_\d(\ov
G,\ov C)$. Such an $A_\d$ exists, as $|P_\d(\ov G,\ov C)|\le
\l<\cont$, while there are $\cont$ subsets of $c_\d$. By Theorem
\typeCo, there is a homogeneous group $H$ with type ${\bf t}$ such
that $A_\d\in P_\d(\ov H,\ov C)$ for every $\d\in S$. If there were a
pure embedding $\phi:H\to G$, then by Theorem \embedding, for almost
every $\d\in S$, $A_\d$ would be in $P_\d(\ov G,\ov C)$. But by the
choice of $A_\d$ this is impossible.\endproof\homNeg

\bigskip
\noindent
{\bf References}

\item{[CK]} C. C Chang and J. J. Keisler {\bf Model Theory}, North
Holland 1973.

\item {[Fu]} L. Fuchs, {\bf Infinite Abelian Groups}, vols. I, II,
Academic press, 1970.

\item{[EM]} P. C. Eklof and A. H. Mekler, {\bf Almost free Modules:  

Set theoretic
methods} North Holland Math. Library, 1990

\item{[Sh-g]} Saharon Shelah, {\bf Cardinal Arithmetic}, to appear in
OUP.

\item{[KjSh 409]} M. Kojman and S. Shelah, 

Existence of universal models, accepted to the {\bf
Journal of Symbolic Logic}

\item {[KjSh 447]}  M. Kojman and S. Shelah, The universality
spectrum of stable unsuperstable theories, accepted to {\bf Annals of
Pure and Applied Logic}.

\item{[GrSh 174]}
R. Grossberg and S. Shelah, On universal locally finite groups,
{\bf Israel J. of Math.},  44, (1983), 289-302.

\bigskip
Department of Mathematics
Carnegie Mellon University
Pittsburgh, PA 15 213
USA

\bigskip
Institute for Mathematics
The Hebrew University of Jerusalem,
Jerusalem 91904
\end